\documentclass[a4paper,12pt]{article}
\usepackage[T2A]{fontenc}
\usepackage[utf8]{inputenc}
\usepackage{amsmath, amsfonts, amssymb, amsthm}
\usepackage[english,russian]{babel}
\usepackage{cmap}
\usepackage{hyperref}
\usepackage{tikz,float}

\newcounter{num}[section]

\newenvironment{statement}
{\refstepcounter{num}%
	\bigskip\noindent\nopagebreak[4]{\bf Утверждение~\arabic{section}.\arabic{num}. }\it}

\newenvironment{consequence}
{\refstepcounter{num}%
	\bigskip\noindent\nopagebreak[4]{\bf Следствие~\arabic{section}.\arabic{num}. }}

\newcommand{\RN}[1]{%
	\textup{\uppercase\expandafter{\romannumeral#1}}%
}

\tikzstyle{vertex}=[circle,fill=white!25,minimum size=20pt,inner sep=0pt]
\tikzstyle{edge} = [draw,thick,-]

\begin{document}
	
	\title{Об уравнениях над полугруппами Брандта}
	\author{М.А. Вахрамеев}
	\maketitle
		
\begin{abstract}
	В работе изучаются уравнения над полугруппой Брандта $B_n$. Для класса уравнений от одной переменной установлено, что  количество несовместных уравнений асимптотически равно $\frac{2}{n^2}$, а среднее число решений таких уравнений асимптотически равно $n^2$.
	
	In this paper, we study equations over Brandt semigroup $B_n$. We compute that the number of unsolvable equations in one variable asymptotically equals $\frac{2}{n^2}$, and the average number of solutions of these equations asymptotically equals $n^2$.
	
\end{abstract}

\section{Введение}

Одним из основных понятий универсальной алгебраической геометрии, развитой в работах Э.Ю.~Данияровой, А.Г.~Мясникова и В.Н.~Ремесленникова (см., например, монографию [1]), является понятие алгебраического множества. Напомним, что алгебраическим множеством над алгебраической системой $A$ называется множество $Y$, для которого существует система уравнений над $A$ с решением $Y$. Среди всех алгебраических множеств фиксированной размерности важную роль играют множества, задаваемые одним уравнением (далее такие множества будем называть 1-алгебраическими). 

Во многих прикладных задачах возникает необходимость порождения случайного 1-алгебраического множества. Наиболее естественное решение данной проблемы заключается в выборе пары термов $t(x),s(x)\in T(x)$ (относительно некоторой вероятностной меры $\mu$ на множестве термов $T(x)$), которая и будет задавать 1-алгебраическое множество $V_A(t(x)=s(x))$. 

В работе [2], где  изучались уравнения над полурешетками, было замечено, что равномерность меры $\mu$ дает далеко не равномерное распределение на классе 1-алгебраических множеств. Данный эффект обусловлен тем, что для 1-алгебраических множеств числа $|Eq(Y)|=|\{t(x)=s(x)\ |\ V_A(t(x)=s(x))=Y\}|$ существенно отличаются друг от друга при разных $Y$. 

Другой способ порождения случайных 1-алгебраических множеств, отчасти исправляющий недостатки предыдущего метода, заключается в следующем. Введем отношение эквивалентности на множестве термов над алгебраической системой $A$:

$$t\sim s \iff t(x)=s(x)\ \text{для любой точки}\ x \in A.$$

Оставляя по одному представителю из каждого класса термов, мы получим множество термов $\bar{T}(x)$, и тогда порождение случайного 1-алгебраического множества заключается в выборе пары термов $t(x),s(x) \in \bar{T}(x)$. 

В настоящей работе второй способ порождения 1-алгебраических множеств применен к алгебраическим множествам, задаваемым уравнениями от одной переменной над полугруппой брандта $B_n$. Уравнение определяется как равенство двух термальных функций. Доказывается, что число классов эквивалентности термальных функций над полугруппой Брандта $B_n$ равно $n^4 + 3n^2 + 3$. Используя данный факт, нам удалось вычислить количество совместных и несовместных уравнений, а также установить, что среднее число решений для уравнений над полугруппой Брандта $B_n$ асимптотически равно $n^2$.

\section{Основные определения}
{\it Полугруппой Брандта} называется полугруппа  $B_n = \{0\} \cup \{(i, j), 1 \leq i,j \leq n \}$, $|B_n| = n^2 + 1$, операция умножения на которой задается следующим образом:

\begin{enumerate}
	\item для любого $b \in B_n$: $b \cdot 0 = 0 \cdot b = 0$;
	\item для любых $(e,k),\ (l,m) \in B_n \setminus \{0\}$:
		\begin{equation}
			(e,k)\cdot(l,m) = 
			\begin{cases}
				(e,m), & \text{если } k = l; \\
				0, & \text{иначе}.
			\end{cases}
		\end{equation}
\end{enumerate}

$B_n$-{\it термом} от переменной $x$ или просто {\it термом} будем называть выражение $w(x) = w_1 \cdot w_2 \dots w_s$ такое, что для любого $1 \leq i \leq s$ либо $w_i \in B_n$, либо $w_i$ есть переменная $x$. Элементы полугруппы $B_n$, входящие в запись $B_n$-терма, мы будем называть {\it константами}. Два терма $s$ и $t$ будем называть равными, если совпадают их записи. Также определим отношение эквивалентности на множестве всех термов : $$s(x) \sim t(x) \iff s(x) = t(x)\ \text{при любом значении $x$.}$$

Функция $f: B_n \rightarrow B_n$ будет называться {\it термальной}, если она задается некоторым $B_n$-термом. Через $T$ обозначим множество всех различных термальных функций над $B_n$. Будем говорить, что термальная функция $f$, заданная термом $w$, эквивалентна терму $v$, если $w \sim v$.

{\it Уравнением} над $B_n$ будем называть равенство двух термальных функций $f_1 = f_2,\ f_1,f_2 \in T$. {\it Решением} уравнения $f_1 = f_2$ называется такое значение переменной $x$, подстановка которого в уравенение обращает его в верное равенство. Уравнение, имеющее хотя бы одно решение, называется {\it совместным }. Уравнение, не имеющее решений, называется {\it несовместным }. Множество всех решений уравнения $f_1 = f_2$ будем обозначать, как $V_{B_n}(f_1 = f_2)$. Множество всех уравнений над $B_n$ будем обозначать, как $Eq_n$.
\section{Эквивалентность термальных функций}

Определим следующие множества термов:

\begin{itemize}
	\item $C_1 = \{ b \in B_n \}$, $|C_1| = n^2 + 1$;
	\item $C_2 = \{x\}$, $|C_2| = 1$;
	\item $C_3 = \{x^2\}$, $|C_3| = 1$;
	\item $C_4 = \{bx,\ b \in B_n \setminus \{0\}\}$,\ $|C_4| = n^2$;
	\item $C_5 = \{xb,\ b \in B_n \setminus \{0\}\}$,\ $|C_5| = n^2$;
	\item $C_6 = \{bxd,\ b,d \in B_n \setminus 0\}$,\ $|C_6| = n^4$.
\end{itemize}

Пусть  $C = C_1 \cup C_2 \cup C_3 \cup C_4 \cup C_5 \cup C_6$. Поскольку множества $C_i$ попарно не пересекаются, то $|C| = |C_1| + |C_2| + |C_3| + |C_4| + |C_5| + |C_6| = n^4 + 3n^2 + 3$.

\begin{statement}
	Для каждой пары термов $t,s \in C$ таких, что $t \sim s$, следует, что $t = s$.
	
	\begin{proof}
	Для того, чтобы убедиться в правильности данного утверждения, достаточно сравнить множества значений для всех термов из $C$:
		
		\begin{itemize}
			\item $t \in C_1,\ t = b,\ b \in B_n$
			\item $t \in C_2,\ t = x$
			\item
				$t \in C_3,\ t = x^2 = 
					\begin{cases}
						(i,i), & \text{если } x = (i,i); \\
						0, & \text{иначе};
					\end{cases}$
			\item $t \in C_4,\ t = bx = (b^1, b^2)x =
					\begin{cases}
						(b^1, i), & \text{если } x = (b^2, i); \\
						0, & \text{иначе};
					 \end{cases}$
			\item $t \in C_5,\ t = xb = x(b^1, b^2) =
					\begin{cases}
						(i, b^2), & \text{если } x = (i, b^1); \\
						0, & \text{иначе};
					\end{cases}$
			\item $t \in C_6,\ t = bxd = (b^1, b^2)x(d^1, d^2) =
					\begin{cases}
						(b^1, d^2), & \text{если } x = (b^2, d^1); \\
						0, & \text{иначе}.
					\end{cases}$
		\end{itemize}
		
	Таким образом, множество $C$ не содержит различных экивалентных термов.
		
	\end{proof}
	
\end{statement}

Далее любая термальная функция $f: B_n \to B_n$ будет отождествляться с соответствующим термом множества $C$.

\begin{statement}
	Для любой термальной функции $f: B_n \to B_n$ существует эквивалентный ей терм $t\in C$.
	
	\begin{proof}
		Пусть функция $f$ задана некоторым термом $w = w_1 \cdot w_2 \dots w_k$. Заменив произведение соседних констант на его результат, терм $w = w_1 \cdot w_2 \dots w_k$ можно привести либо к виду $v = v_1 \cdot v_2 \dots v_s$, не содержащему в своей записи двух подряд идущих констант, либо к 0. Если $f \sim 0$, то утверждение автоматически доказано, поскольку $0 \in C_1$. Рассмотрим полученный выше терм $v = v_1 \cdot v_2 \dots v_s$.
		
		Если $s \leq 2$, то $v$ - это терм одного из следующих видов: $b$, $x$, $x^2$, $bx$, $xb$. Все такие термы содержатся во множестве $C$, cледовательно, для случая $s \leq 2$ утверждение также верно.
		
		Теперь рассмотрим случай, когда $s > 2$. Заметим, что \newline
		$$x^s = 
			\begin{cases}
				(i,i), & \text{если } x = (i,i); \\
				0, & \text{иначе},
			\end{cases}
		$$
		следовательно, $x^s \sim x^2$. Таким образом, если $v = x^s$, то утвержедение доказано. Поэтому предположим, что $v$ содержит хотя бы одну константу. Докажем, что любой такой терм либо тождественно равен 0, либо не равен нулю только на одной точке. Заметим, что $v$ в своей записи обязательно будет содержать одно из следующих выражений:
		
		\begin{enumerate}
			\item $(b^1,b^2) \cdot x \cdot x = 
				\begin{cases}
					(b^1, b^2), & \text{если } x = (b^2, b^2); \\
					0, & \text{иначе};
				\end{cases}$	
			\item $x \cdot (b^1,b^2) \cdot x = 
				\begin{cases}
					(b^2, b^1), & \text{если } x = (b^2, b^1); \\
					0, & \text{иначе};
				\end{cases}$	
			\item $x \cdot x \cdot (b^1,b^2) = 
				\begin{cases}
					(b^1, b^2), & \text{если } x = (b^1, b^1); \\
					0, & \text{иначе};
				\end{cases}$	
			\item $(b^1,b^2) \cdot x \cdot (d^1, d^2) = 
				\begin{cases}
					(b^1, d^2), & \text{если } x = (b^2, d^1); \\
					0, & \text{иначе}.
				\end{cases}$
		\end{enumerate}
		
		Любое из произведений 1-4 не равняется нулю лишь при единственном значении $x \neq 0$. Следовательно, терм $v$ либо тождественно равен 0, либо не равен 0 максимум при одном значении $x \neq 0$. Если терм $v$ отличен от нуля только при $x = (a^1, a^2)$ и $v((a^1, a^2)) = (e^1, e^2)$, то $v \sim (e^1, a^1)x(a^2, e^2) \in C$.		
	\end{proof}
	
\end{statement}

\begin{consequence}
	Число различных термальных функций над $B_n$ равно $n^4 + 3n^2 + 3$.
\end{consequence}

\section{Среднее число решений уравнений над $B_n$}	
	
Согласно определению, уравнением над $B_n$ называется равенство двух термальных функций $\ f_1,f_2 \in T$: $f_1(x) = f_2(x)$. Так как $|T| = n^4 + 3n^2 + 3$, то $|Eq_n| = (n^4 + 3n^2 + 3)^2 = n^8 + 6n^6 + 15n^4 + 18n^2 + 9$.

Разобъем множество $Eq_n$ на подмножества, каждое из которых будет состоять из уравнений определенного типа. Тип уравнения будет определяться принадлежностью термальных функций, его задающих, к множествам $C_1, \dots , C_6$. Таким образом, каждое уравнение из $Eq_n$ относится одному из $6\cdot6 = 36$ типов. Найдем количество решений для уравений каждого типа.

\subsection{$f_1 \sim t_1 \in C_1, f_2 \sim t_2 \in C_1:\ b_1 = b_2,\ b_1,b_2 \in B_n$}
Число уравнений данного типа равно $|C_1|^2 = (n^2 + 1)^2 = n^4 + 2n^2 + 1$. Количество решений для уравнения, левая и правая часть которого задаются одним и тем же термом из $C_1$, равно числу элементов в $B_n$ т.е. $n^2 + 1$. Таких уравнений ровно $|C_1| = n^2 + 1$ штук. Для остальных $n^4 + n^2$ уравнений решений нет.

\begin{table}[H]
	\begin{tabular}{|c|c|c|}
		\hline
		Кол-во решений\ \ \ \ & 0 & $n^2 + 1$ \\
		\hline
		Кол-во уравнений & $n^4 + n^2$ & $n^2 + 1$ \\
		\hline
	\end{tabular}
\end{table}
	
\subsection{$\ f_1 \sim t_1 \in C_1,\ f_2 \sim t_2 \in C_2:\ b = x,\ b \in B_n$}

Число уравнений данного типа равно $|C_1|\cdot|C_2| = n^2 + 1$. Каждое уравнение имеет единственное решение. Учтем также уравнения симметричного вида, где $f_1 \sim t_1 \in C_2$,\ $f_2 \sim t_2 \in C_1$.

\begin{table}[H]
	\begin{tabular}{|c|c|}
		\hline
		Кол-во решений\ \ \ \  & 1 \\
		\hline
		Кол-во уравнений & $2(n^2 + 1)$ \\
		\hline
	\end{tabular}
\end{table}

\subsection{$f_1 \sim t_1 \in C_1,\ f_2 \sim t_2 \in C_3:\ b = x^2,\ b \in B_n$}

Число уравнений данного типа равно $|C_1|\cdot|C_3| = n^2 + 1$.
Так как
$$
x^2 =
\begin{cases}
(t, t), & \text{если } x = (t, t),\ 1 \leq t \leq n; \\
0, & \text{иначе};
\end{cases}
$$
то:
\begin{enumerate}
	\item[-] если $b \neq 0$\ и $b \neq (i,i),\ 1 \leq i \leq n$ ($n^2 + 1 - 1 - n = n^2 - n$ уравнений), то решений нет;
	\item[-] если $b = (i,i), \ 1 \leq i \leq n$ ($n$ уравнений), то существует ровно 1 решение - $(i,i)$;
	\item[-] если $b = 0$ (1 уравнение), то решением будет любое значение $x$, кроме $x = (i,i),\ 1 \leq i \leq n$, то есть будет существовать $n^2 + 1 - n$ решений.
\end{enumerate}

Учитывая также уравнения симметричного типа $x^2 = b$, получаем:

\begin{table}[H]
	\begin{tabular}{|c|c|c|c|}
		\hline
		Кол-во решений\ \ \ \ & 0 & $1$ & $n^2 - n + 1$\\
		\hline
		Кол-во уравнений & $2(n^2 - n)$ & $2n$ & $2$\\
		\hline
	\end{tabular}
\end{table}

\subsection{$f_1 \sim t_1 \in C_1,\ f_2 \sim t_2 \in C_4:\ b_1 = b_2x,\ b_1, b_2 \in B_n, b_2 \neq 0$}

Количество уравнений данного типа равно $|C_1|\cdot|C_4| = (n^2 + 1)n^2$.

Пусть $b_1 \neq 0$ и $b_1 = (b_1^1, b_1^2)$, а $b_2 = (b_2^1, b_2^2)$, и уравнение запишется, как $(b_1^1, b_1^2) = (b_2^1, b_2^2)x$. Если $b_1^1 = b_2^1$\ ($n \cdot n \cdot n = n^3$ уравнений), то будет существовать только одно решение $x = (b_2^2, b_1^2)$. Если $b_1^1 \neq b_2^1$ ($n \cdot n \cdot (n - 1) \cdot n = n^4 - n^3$ уравнений), то решений нет.

Пусть теперь $b_1 = 0$, $b_2 = (b_2^1, b_2^2)$ и уравнение запишется как $0 = (b_2^1, b_2^2)x$. Тогда любое значение $x$, кроме $x = (b_2^2, i),\ 1 \leq i \leq n$, будет являться решением. Таким образом, для любого из $n \cdot n = n^2$ уравнений будет существовать $n^2 + 1 - n$ решений.

Аналогичные рассуждения верны и для уравнений, где $f_1 \sim t_1 \in C_1,\ f_2 \sim t_2 \in C_5$, то есть для уравнений вида $\ b_1 = xb_2,\ b_1, b_2 \in B_n, b_2 \neq 0$. Также, учитывая уравнения симметричных типов $b_2x = b_1$ и $xb_2 = b_1$, получаем:

\begin{table}[H]
	\begin{tabular}{|c|c|c|c|}
		\hline
		Кол-во решений\ \ \ \ & 0 & $1$ & $n^2 - n + 1$\\
		\hline
		Кол-во уравнений & $4(n^4 - n^3)$ & $4n^3$ & $4n^2$\\
		\hline
	\end{tabular}
\end{table}

\subsection{$f_1 \sim t_1 \in C_1,\ f_2 \sim t_2 \in C_6:\ b = cxd,\ b,c,d \in B_n,\ c,d \neq 0$}

Число уравнений данного типа равно $|C_1|\cdot|C_6| = (n^2 + 1)n^4$.

Пусть $b = (b^1, b^2),\  c = (c^1, c^2),\ d = (d^1, d^2)$. Если $(c^1, d^2) = (b^1, b^2)$ ($n \cdot n \cdot n \cdot n = n^4$ уравнений), то существует ровно одно решение $x = (c^2, d^1)$. В противном случае, если $(c^1, d^2) \neq (b^1, b^2)$ ($n^6 - n^4$ уравнений, где $n^6$ - число вариантов выбрать $b,c,d,\ \text{т.ч.}\ b,c,d \neq 0$, $n^4$ - число уравнений, где $(c^1, d^2) = (b^1, b^2)$), то решений нет.

Теперь рассмотрим случай, когда $b = 0$ и уравнение принимает вид $0 = cxd$. Количество таких уравнений равно $|C_6| = n^4$. Для каждого из таких уравнений любое значение $x$, обращающее $cxd$ в ноль, является решением. Таким образом, любое такое уравнение будет иметь $n^2 + 1 - 1 = n^2$ решений.

Учитывая также уравнения симметричного вида $cxd = b$, получаем:

\begin{table}[H]
	\begin{tabular}{|c|c|c|c|}
		\hline
		Кол-во решений\ \ \ \ & 0 & $1$ & $n^2$\\
		\hline
		Кол-во уравнений & $2(n^6 - n^4)$ & $2n^4$ & $2n^4$\\
		\hline
	\end{tabular}
\end{table}

\subsection{$f_1 \sim t_1 \in C_2,\ f_2 \sim t_2 \in C_2:\ x = x$}

Одно уравнение, $n^2 + 1$ решений.

\begin{table}[H]
	\begin{tabular}{|c|c|}
		\hline
		Кол-во решений\ \ \ \ & $n^2 + 1$\\
		\hline
		Кол-во уравнений & 1\\
		\hline
	\end{tabular}
\end{table}

\subsection{$f_1 \sim t_1 \in C_2,\ f_2 \sim t_2 \in C_3:\ x = x^2$}

Решениями будут элементы $B_n$ вида $(i, i),\ 1 \leq i \leq n$, и ноль. Учтем также уравнение $x^2 = x$.

\begin{table}[H]
	\begin{tabular}{|c|c|}
		\hline
		Кол-во решений\ \ \ \ & $n + 1$\\
		\hline
		Кол-во уравнений & 2\\
		\hline
	\end{tabular}
\end{table}

\subsection{$f_1 \sim t_1 \in C_2,\ f_2 \sim t_2 \in C_4:\ x = bx,\ b \in B_n,\ b \neq 0$}

Количество уравнений данного типа равно $|C_2| \cdot |C_4| = n^2$. Отметим, что для любого из этих уравнений $x = 0$ является решением.

Далее, пусть $b = (i, i),\ 1 \leq i \leq n$, и уравнения имеют вид $x = (i,i)x$ (всего $n$ уравнений). Тогда, кроме 0, решениями будут также $x = (i, z),\ 1 \leq z \leq n$.  Таким образом, каждое из $n$ уравнений будет иметь $n + 1$ решение.

Пусть теперь $b = (i, j),\ 1 \leq i,j \leq n,\ i \neq j$, и уравнения имеют вид $x = (i,j)x$ ($n^2 - n$ уравнений). Других решений, кроме $x = 0$, здесь не будет, поскольку, чтобы правая часть уравнения не равнялась нулю, $x$ должен иметь первую координату равную $j$, но, так как произведение $(i,j)x$ будет иметь первую координату равную $i$, при $j \neq i$ равенство левой и правой частей уравнения невозможно. Следовательно, каждое из $n^2 - n$ уравнений будет иметь единственное решение: $x = 0$.

Заметим, что аналогичные рассуждения верны и для уравнений, где $f_1 \sim t_1 \in C_2,\ f_2 \sim t_2 \in C_5$, то есть для уравнений вида $x = xb,\ b \in B_n,\ b \neq 0$. Учитывая также уравнения симметричных типов $bx = x$ и $xb = x$, получаем:

\begin{table}[H]
	\begin{tabular}{|c|c|c|}
		\hline
		Кол-во решений\ \ \ \ & 1 & $n + 1$\\
		\hline
		Кол-во уравнений & $4(n^2 - n)$ & $4n$\\
		\hline
	\end{tabular}
\end{table}

\subsection{$f_1 \sim t_1 \in C_2,\ f_2 \sim t_2 \in C_6:\ x = bxd,\ b,d \in B_n \setminus 0$}

Количество уравнений данного типа равно $|C_2|\cdot|C_6| = n^4$. Для любого из этих уравнений $x = 0$ является решением.
 
Пусть $b = (b^1, b^2),\ d = (d^1, d^2)$. Очевидно, что, если $x \neq 0$ и $x \neq (b^2, d^1)$, то $x$ не может быть решением. Допустим, $x = (b^2, d^1)$. Тогда, если $(b^2, d^1) = (b^1, d^2)$ ($n^2$ уравнений), то $x = (b^2, d^1)$ - решение, в противном случае (для $n^4 - n^2$ уравнений) - нет. Учитывая уравнения симметричного типа $bxd = x$, получаем:

\begin{table}[H]
	\begin{tabular}{|c|c|c|}
		\hline
		Кол-во решений\ \ \ \ & 1 & 2\\
		\hline
		Кол-во уравнений & $2(n^4 - n^2)$ & $2n^2$\\
		\hline
	\end{tabular}
\end{table}

\subsection{$f_1 \sim t_1 \in C_3,\ f_2 \sim t_2 \in C_3:\ x^2 = x^2$}

Одно уравнение, $n^2 + 1$ решений.

\begin{table}[H]
	\begin{tabular}{|c|c|}
		\hline
		Кол-во решений\ \ \ \ & $n^2 + 1$\\
		\hline
		Кол-во уравнений & 1\\
		\hline
	\end{tabular}
\end{table}

\subsection{$f_1 \sim t_1 \in C_3,\ f_2 \sim t_2 \in C_4:\ x^2 = bx,\ b \in B_n,\ b \neq 0$}

Количество уравнений данного типа равно $|C_3| \cdot |C_4| = n^2$. Для любого такого уравнения $x = 0$ является решением.

Пусть $b = (b^1, b^2)$ и уравнения имеют вид $x^2 = (b^1, b^2)x$. Допустим, что $b^1 = b^2 = b'$ ($n$ уравнений). Тогда любой $x$ вида $(b', t)$,\ $1 \leq t \leq n$, за исключением $(b', b')$, обращает левую часть уравнения в ноль, в то время как правая часть уравнения в ноль не обращается. Следовательно, $x = (b', t),\ 1 \leq t \leq n,\ t \neq b'$ решениями не являются. При этом $x = (b', b')$, очевидно, есть решение. Далее, все $x$ вида $(i, i),\ 1 \leq i \leq n,\ i \neq b'$ обращают в ноль правую часть уравнения, при этом не зануляя левую, поэтому решениями не будут. Наконец, все остальные значения $x \neq (i,i),\ x \neq (b', t)$ одновременно зануляют обе части уравнения, следовательно являются решениями. Таким образом, каждое из $n$ уравнений будет иметь $n^2 + 1 - (n - 1) - (n - 1) = n^2 - 2n + 3$ решений.

Теперь рассмотрим случай, когда $b^1 \neq b^2$ ($n^2 - n$ уравнений). Если $x$ имеет вид $(b^2, z),\ 1 \leq z \leq n,\ z \neq b^2$, то левая часть уравнения зануляется, а правая нет, следовательно, данные $x$ решениями не будут. При $x = (b^2, b^2)$ обе части уравнения ненулевые, но при этом не равны друг другу, а значит $x = (b^2, b^2)$ не является решением. Если $x$ имеет вид $(i, i),\ 1 \leq i \leq n,\ i \neq b^2$, то левая часть уравнения ненулевая, а правая обращается в ноль, значит, данные $x$ также не могут быть решениями. Все остальные значения $x \neq (i, i),\ x \neq (b^2, z)$ обращают в ноль обе части уравнения, поэтому являются решениями. Получаем, что любое из $n^2 - n$ уравнений будет иметь $n^2 + 1 - (n - 1) - 1 - (n - 1) = n^2 - 2n + 2$ решений.

Аналогичные рассуждения верны и для уравнений, где $f_1 \sim t_1 \in C_3,\ f_2 \sim t_2 \in C_5$, то есть для уравнений вида $x^2 = xb,\ b \in B_n,\ b \neq 0$. Учитывая также уравнения симметричных типов $bx = x^2$ и $xb = x^2$, получаем:

\begin{table}[H]
	\begin{tabular}{|c|c|c|}
		\hline
		Кол-во решений\ \ \ \ & $n^2 - 2n + 3$ & $n^2 - 2n + 2$\\
		\hline
		Кол-во уравнений & $4n$ & $4(n^2 - n)$\\
		\hline
	\end{tabular}
\end{table}

\subsection{$f_1 \sim t_1 \in C_3,\ f_2 \sim t_2 \in C_6:\ x^2 = bxd,\ b,d \in B_n \setminus 0$}

Количество уравнений данного типа равно $|C_3|\cdot|C_6| = n^4$. Для любого их этих уравнений $x = 0$ является решением.

Пусть $b = (b^1, b^2),\ d = (d^1, d^2)$. Для начала рассмотрим случай, когда $b^2 = d^1 = z$ ($n \cdot n^2 = n^3$ уравнений). 

Допустим, что $b^1 = d^2 = z$ ($n$ уравнений). Нетрудно видеть, что  $x = (z, z)$ есть решение. Далее, если $x = (i, i) \in B_n,\ i \neq z$, то правая часть уравнения зануляется, а левая нет, поэтому данные значения $x$ решениями не будут. Наконец, любое $x \neq (i, i) \in B_n$ зануляет обе части уравнения, следовательно, решением является. Таким образом, каждое  из $n$ уравнений будет иметь $n^2 + 1 - (n - 1) = n^2 - n + 2$ решений.

Пусть теперь условие $b^1 = d^2 = z$ не выполнено ($n^3 - n$ уравнений). Тогда, если $x = (i, i) \in B_n, i \neq z$, то правая часть уравнения зануляется, а левая нет, поэтому такие $x$ не будут решениями. В случае $x = (z, z)$ обе части уравнения не равны нулю, но при этом не совпадают, следовательно $x = (z,z)$ также не решение. Все остальные $x \neq (i,i)$ решениями являются, поскольку зануляют обе части уравнения. Итого, каждое из $n^3 - n$ уравнений будет иметь $n^2 + 1 - n$ решений.

Теперь рассмотрим случай, когда $b^2 \neq d^1$ ($n^4 - n^3$ уравнений). При $x = (i,i)\ \in B_n$ правая часть уравнения равна нулю, левая нулю не равна, поэтому данные значения $x$ решениями не будут. Если $x = (b^2, d^1)$, то, наоборот, левая часть уравнения занулится, а правая нет, поэтому $x = (b^2, d^1)$ также не решение. Остальные значения $x$ являются решениями, поскольку обращают в ноль обе части уравнения. Таким образом, для любого из $n^4 - n^3$ уравнений будет существовать $n^2 + 1 - n - 1 = n^2 - n$ решений.

Учитывая уравнения симметричного типа $bxd = x^2$, получаем:

\begin{table}[H]
	\begin{tabular}{|c|c|c|c|}
		\hline
		Кол-во решений\ \ \ \ & $n^2 - n$ & $n^2 - n + 1$ & $n^2 - n + 2$\\
		\hline
		Кол-во уравнений & $2(n^4 - n^3)$ & $2(n^3 - n)$ & $2n$\\
		\hline
	\end{tabular}
\end{table}

\subsection{$f_1 \sim t_1 \in C_4,\ f_2 \sim t_2 \in C_4:\ b_1x = b_2x,\ b_1,b_2 \in B_n,\ b_1 \neq 0,b_2 \neq 0$}

Число уравнений данного типа равно $|C_4|\cdot|C_4| = n^2 \cdot n^2 = n^4$. Для любого такого уравнения $x = 0$ является решением.

Если $b_1 = b_2$ ($n^2$ уравнений), то любое $x \in B_n$ является решением. Таким образом, получаем, что каждое из $n^2$ уравнений будет иметь $n^2 + 1$ решение. 

Пусть теперь $b_1 \neq b_2$, $b_1 = (b_1^1, b_1^2),\ b_2 = (b_2^1, b_2^2)$, и уравнения имеют вид $(b_1^1, b_1^2)x = (b_2^1, b_2^2)x$. Если $x \neq (b_1^2, y)$ и $x \neq (b_2^2, z),\ 1 \leq y,z \leq n$, то $x$ зануляет левую и правую части уравнения, следовательно явялется решением. 

Далее, предположим, что $b_1^2 \neq b_2^2$ ($n \cdot n \cdot n \cdot (n - 1) = n^4 - n^3$ уравнений). Тогда все $x = (b_1^2, y)$ зануляют правую часть уравнения и не зануляют левую, поэтому решениями быть не могут. Наоборот, $x = (b_2^2, z)$ зануляют левую часть уравнения, не зануляют правую, следовательно, также не являются решениями. Итого, для каждого из $n^4 - n^3$ уравнений будет существовать ровно $n^2 + 1 - n - n = n^2 - 2n + 1$ решение.

Теперь рассмотрим случай $b_1^2 = b_2^2 = b'$ ($n^4 - n^2 - (n^4 - n^3) = n^3 - n^2$ уравнений). При $x = (b', t),\ 1 \leq t \leq n$ и левая и правая части уравнения нулю не равны, но при этом они и не совпадают, так как $b_1^1 \neq b_2^1$. Следовательно, $x = (b', t)$ решениями не являются. Получаем, что каждое из $n^3 - n^2$ уравнений будет иметь $n^2 + 1 - n$ решений.

Аналогичные рассуждения верны и для уравнений вида $xb_1 = xb_2,\ b_1,b_2 \in B_n,\ b_1 \neq 0,b_2 \neq 0$. Таким образом, получаем:

\begin{table}[H]
	\begin{tabular}{|c|c|c|c|}
		\hline
		Кол-во решений\ \ \ \ & $n^2 +1$ & $n^2 - n + 1$ & $n^2 - 2n + 1$\\
		\hline
		Кол-во уравнений & $2n^2$ & $2(n^3 - n^2)$ & $2(n^4 - n^3)$\\
		\hline
	\end{tabular}
\end{table}

\subsection{$f_1 \sim t_1 \in C_4,\ f_2 \sim t_2 \in C_5:\ b_1x = xb_2,\ b_1,b_2 \in B_n,\ b_1 \neq 0,b_2 \neq 0$}

Количество уравнений данного типа равно $|C_4|\cdot|C_5| = n^2 \cdot n^2 = n^4$. Для любого из этих уравнений $x = 0$ является решением.

Пусть $b_1 = (b_1^1, b_1^2),\ b_2 = (b_2^1, b_2^2)$, и уравнения имеют вид $(b_1^1, b_1^2)x = x(b_2^1, b_2^2)$. Все значения $x \neq (b_1^2, t)$ и $x \neq (z, b_2^1)$, $1 \leq z,t \leq n$, обращают обе части уравнения в ноль, поэтому являются решениями. При этом, $x =(b_1^2, t),\ t \neq b_2^1$, зануляют правую часть и не зануляют левую, следовательно решениями не будут. Значения $x = (z, b_2^1),\ z \neq b_1^2$, напротив, обращают в ноль левую часть, оставляя ненулевой правую, поэтому также не являются решениями. Осталось рассмотреть случай, когда $x = (b_1^2, b_2^1)$. Подставим данное значение $x$ в уравнение: 
	$$(b_1^1, b_1^2)x = x(b_2^1, b_2^2)$$ 
	$$(b_1^1, b_1^2)(b_1^2, b_2^1) = (b_1^2, b_2^1)(b_2^1, b_2^2)$$
	$$(b_1^1,b_2^1) = (b_1^2, b_2^2).$$

Если $b_1^1 = b_1^2$ и $b_2^1 = b_2^2$ ($n \cdot n = n^2$ уравнений), то $x = (b_1^2, b_2^1)$ - решение, иначе ($n^4 - n^2$ уравнений) - нет.

Таким образом, каждое из $n^2$ уравнений будут иметь $n^2 + 1 - (n - 1) - (n - 1) = n^2 - 2n + 3$ решения, а каждое из $n^4 - n^2$ уравнений будет иметь $n^2 + 1 - (n - 1) - (n - 1) - 1 = n^2 - 2n + 2$ решения. Учитывая уравнения симметричного типа $xb_2 = b_1x$, получаем:

\begin{table}[H]
	\begin{tabular}{|c|c|c|}
		\hline
		Кол-во решений\ \ \ \ & $n^2 - 2n + 3$ & $n^2 - 2n + 2$\\
		\hline
		Кол-во уравнений & $2n^2$ & $2(n^4 - n^2)$\\
		\hline
	\end{tabular}
\end{table}

\subsection{$f_1 \sim t_1 \in C_4,\ f_2 \sim t_2 \in C_6:\ bx = cxd,\ b,c,d \in B_n \setminus 0$}

Количество уравнений данного вида равно $|C_4|\cdot|C_6| = n^2 \cdot n^4 = n^6$. Для любого уравнения данного вида $x = 0$ является решением.

Пусть $b = (b^1, b^2),\ c = (c^1, c^2),\ d = (d^1, d^2)$, и уравнения имеют вид $(b^1, b^2)x = (c^1, c^2)x(d^1, d^2)$. Все значения $x \neq (b^2, t),\ 1 \leq t \leq n,\ x \neq (c^2, d^1)$ обращают обе части уравнения в ноль, поэтому являются решениями.

Допустим, что $c^2 = b^2$ ($n \cdot n \cdot n \cdot n \cdot n = n^5$ уравнений). Тогда $x = (b^2, t),\ t \neq d^1$ решениями не будут, так как зануляют правую часть уравнения и не зануляют левую. Осталось рассмотреть случай, когда $x = (b^2, d^1)$. Если $ (c^1, d^2) = (b^1, d^1)$ ($n \cdot n \cdot n = n^3$ уравнений), то $x = (b^2, d^1)$ - решение. В противном случае ($n^5 - n ^3$ уравнений), $x = (b^2, d^1)$ решением не является. Таким образом, получаем, что каждое из $n^3$ уравнений будет иметь $n^2 + 1 - (n - 1) = n^2 -n + 2$ решения, а каждое из $n^5 - n^3$ уравнений будет иметь $n^2 + 1 - n$ решений.

Пусть теперь $c^2 \neq b^2$ ($n^6 - n^5$ уравнений). Левая и правая части уравнения обращаются в ноль при всех значениях $x$, кроме $x = (b^2, t),\ 1 \leq t \leq n$ и $x = (c^2, d^1)$. При этом $x = (b^2, t)$ и $x = (c^2, d^1)$ решениями не являются, так как зануляют одну из частей уравнения, оставляя ненулевой другую. Следовательно, для каждого из $n^6 - n^5$ уравнений будет существовать $n^2 + 1 - n - 1 = n^2 - n$ решений.

Аналогичные рассуждения верны и для уравнений, где $f_1 \sim t_1 \in C_5,\ f_2 \sim t_2 \in C_6$, то есть уравнений вида $xb = cxd,\ b,c,d \in B_n \setminus 0$. Учитывая также уравнения симметричных типов $cxd = bx$ и $cxd = xb$, получаем:

\begin{table}[H]
	\begin{tabular}{|c|c|c|c|}
		\hline
		Кол-во решений\ \ \ \ & $n^2 - n$ & $n^2 - n + 1$ & $n^2 - n + 2$\\
		\hline
		Кол-во уравнений & $4(n^6 - n^5)$ & $4(n^5 - n^3)$ & $4n^3$\\
		\hline
	\end{tabular}
\end{table}

\subsection{$f_1 \sim t_1 \in C_6,\ f_2 \sim t_2 \in C_6:\ bxc = dxe,\ b,c,d,e \in B_n \setminus 0$}

Количество уравнений данного типа равно $|C_6|\cdot|C_6| = n^4 \cdot n^4 = n^8$. Для любого такого уравнения $x = 0$ является решением.

Функции $f_1$ и $f_2$ равны всюду, за исключеним, быть может, двух точек. Это означает, что уравнение гарантированно будет иметь $n^2 + 1 - 2 = n^2 - 1$ решений. Дополнительно будет существовать еще 2 решения (суммарно $n^2 +1$) только в случае, если $f_1$ и $f_2$ полностью совпадают ($n^4$ уравнений). Наконец, если функции $f_1$ и $f_2$ не равны нулю на одной точке, но при этом не совпадают, уравнение будет иметь $n^2 + 1 - 1 = n^2$ решений. Количество таких уравнений будет равно $n^2 \cdot n^2 \cdot (n^2 - 1) = n^6 - n^4$ (здесь $n^2$ число вариантов выбрать общую точку, не переходящую в ноль, умножается на количество возможных образов этой точки для одной из функций, умножается на количество возможных образов общей точки для другой функции, с учетом того, чтобы функции не совпали). Таким образом, получаем:

\begin{table}[H]
	\begin{tabular}{|c|c|c|c|}
		\hline
		Кол-во решений\ \ \ \ & $n^2 - 1$ & $n^2$ & $n^2 + 1$\\
		\hline
		Кол-во уравнений & $n^8 - n^4 - (n^6 - n^4) = n^8 - n^6$ & $n^6 - n^4$ & $n^4$\\
		\hline
	\end{tabular}
\end{table}

\bigskip

Мы рассмотрели все возможные типы уравнений. На основе полученных результатов вычислим среднее число решений случайно выбранного уравнения из $Eq_n$:

\begin{table}[H]
	\begin{tabular}{|c|c|}
		\hline
		Количество решений & Количество уравнений \\
		\hline
		0 & $2n^6 + 3n^4 - 4n^3 + 3n^2 -2n$ \\
		\hline
		1 & $4n^4 + 4n^3 + 4n^2 -2n + 2$ \\
		\hline
		2 & $2n^2$ \\
		\hline
		$n + 1$ & $4n + 2$ \\
		\hline
		$n^2 - 2n + 1$ & $2n^4 - 2n^3$ \\
		\hline
		$n^2 - 2n + 2$ & $2n^4 + 2n^2 - 4n$ \\
		\hline
		$n^2 - 2n + 3$ & $2n^2 + 4n$ \\
		\hline
		$n^2 - n$ & $4n^6 - 4n^5 + 2n^4 - 2n^3$ \\
		\hline
		$n^2 - n + 1$ & $4n^5 + 2n^2 -2n + 2$ \\
		\hline
		$n^2 - n + 2$ & $4n^3 + 2n$ \\
		\hline
		$n^2 - 1$ & $n^8 - n^6$ \\
		\hline
		$n^2$ & $n^6 + n^4$ \\
		\hline
		$n^2 + 1$ & $n^4 + 3n^2 + 3$ \\
		\hline
	\end{tabular}
\end{table} 

$$\overline{Sol}(B_n) = \frac{1}{n^8 + 6n^6 + 15n^4 + 18n^2 + 9} \cdot (4n^4 + 4n^3 + 4n^2 - 2n + 2 + 2\cdot2n^2 + (n + 1)(4n + 2) + $$
$$+ (n^2 - 2n + 1)(2n^4 - 2n^3) + (n^2 - 2n + 2)(2n^4 + 2n^2 - 4n) + (n^2 - 2n + 3)(2n^2 + 4n) + $$
$$+ (n^2 - n)(4n^6 - 4n^5 + 2n^4 - 2n^3) + (n^2 - n +1)(4n^5 + 2n^2 - 2n + 2) + (n^2 - n + 2)(4n^3 + 2n) +  $$
$$+ (n^2 - 1)(n^8 - n^6) + n^2(n^6 + n^4) + (n^2 + 1)(n^4 + 3n^2 + 3)) = $$
$$= \frac{n^{10} + 3n^8 - 4n^7 + 9n^6 - 6n^5 + 22n^4 + 32n^2 + 8n + 9}{n^8 + 6n^6 + 15n^4 + 18n^2 + 9} \sim n^2.$$

\bigskip

Также из таблицы видно, что доля несовместных уравнений равна $\frac{2n^6 + 3n^4 - 4n^3 + 3n^2 -2n}{n^8 + 6n^6 + 15n^4 + 18n^2 + 9} \sim \frac{2}{n^2}$.

\end{document}